\documentclass[12pt, reqno]{amsart}

\usepackage{amssymb}
\usepackage{amsmath}

\usepackage{hyperref}
\usepackage{color}

\def\Z{\mathbb{Z}}
\def\F{\mathcal{F}}

\begin{document}

\title{Small examples of mosaics of combinatorial designs}

\author{Vedran Kr\v{c}adinac}

\address{Faculty of Science, University of Zagreb, Bijeni\v{c}ka cesta~$30$, HR-$10000$ Zagreb, Croatia}

\email{vedran.krcadinac@math.hr}

\thanks{This work has been supported by the Croatian Science Foundation
under the project $9752$.}



\date{May 21, 2024}

\begin{abstract}
We give the first example of a mosaic of three combinatorial designs
with distinct parameters $2$-$(13,3,1)$, $2$-$(13,4,2)$, and
$2$-$(13,6,5)$. Furthermore, we give examples of mosaics of
$2$-$(9,3,2)$ designs that are not resolvable, thereby answering
a question posed by M.~Wiese and H.~Boche. Finally, we give an
example of a mosaic of projective planes of order~$3$ that cannot
be obtained by tiling groups with difference sets.
\end{abstract}

\maketitle

\section{Introduction}\label{sec1}

We use standard notation and facts from design theory~\cite{BJL86, CD07}.
Let $t_i$-$(v,k_i,\lambda_i)$, $i=1,\ldots,c$ be parameters of combinatorial designs
with the same numbers~$v$ of points and $b$ of blocks, and $k_1+\ldots+k_c=v$.
A \emph{mosaic of combinatorial designs} with parameters
$$t_1\mbox{-}(v,k_1,\lambda_1) \oplus \cdots \oplus t_c\mbox{-}(v,k_c,\lambda_c)$$
is a $v\times b$ matrix with entries from $\{1,\ldots,c\}$ such that the $i$-entries
represent incidences of a $t_i$-$(v,k_i,\lambda_i)$ design. Here, $c$ is the number
of \emph{colors} of the mosaic. A mosaic can be thought of as a system of~$c$ designs
with points $\{1,\ldots,v\}$, blocks $\{1,\ldots,b\}$, and incidence relations being
disjoint and covering the whole Cartesian product $\{1,\ldots,v\}\times \{1,\ldots,b\}$.
Mosaics of combinatorial designs were introduced in~\cite{GGP18} and have applications
for information-theoretic security~\cite{WB22, WB24}. Special mosaics of symmetric
designs arise from tilings of groups with difference sets~\cite{CKZ15}. These objects
have been used to construct asynchronous channel hopping systems for cognitive
radio networks~\cite{CZ17, GLZ19}.

Since complementing the blocks of a $t$-$(v,k,\lambda)$ design gives a
$t$-$(v,v-k,\overline{\lambda})$ design for
$\overline{\lambda}=\lambda {v-t\choose k} / {v-t\choose k-t}$, mosaics
with $c=2$ colors are essentially just designs. To avoid trivialities,
we assume $c\ge 3$ and that the strength of each design is $t\ge 2$.
In~\cite{GGP18}, a construction of a mosaic with $c=v/k$ colors from
a resolvable $t$-$(v,k,\lambda)$ design is given, comprising~$c$ isomorphic
copies of the same design. We refer to mosaics with identical parameters
$$t\mbox{-}(v,k,\lambda) \oplus \cdots \oplus t\mbox{-}(v,k,\lambda)$$
as \emph{homogenous}. The authors of~\cite{GGP18} ask for an
inhomogenous mosaic of symmetric designs with parameters
$$2\mbox{-}(31,6,1)  \oplus 2\mbox{-}(31,10,3)\oplus 2\mbox{-}(31,15,7).$$
We do not know if such a mosaic exists, but in Section~\ref{sec2} we
provide the first nontrivial example of an inhomogenous mosaic for
$$2\mbox{-}(13,3,1)\oplus 2\mbox{-}(13,4,2)\oplus 2\mbox{-}(13,6,5).$$

In~\cite{WB24}, the authors ask whether the construction
from~\cite{GGP18} based on resolvable designs is the only way to
obtain homogenous mosaics. We answer the question in the negative
in Section~\ref{sec3} by giving examples containing designs that
are not resolvable.

Strictly speaking, mosaics of symmetric designs are never homogenous,
because $k$ cannot divide $v$. This follows from the necessary existence
condition $\lambda(v-1)=k(k-1)$. However, if $k$ divides $v-1$, a mosaic
of symmetric designs with parameters
$$2\mbox{-}(v,k,\lambda) \oplus \cdots \oplus 2\mbox{-}(v,k,\lambda)\oplus 2\mbox{-}(v,1,0)$$
may be possible. Such mosaics are obtained by developing tilings
of groups with difference sets and will be considered homogenous mosaics
of symmetric designs with $c=(v-1)/k$ colors. Incidences of the trivial
$2$-$(v,1,0)$ design are written as zeros on the diagonal of the
$v\times v$ matrix. In~\cite{CKZ15}, tilings of groups with $(7,3,1)$,
$(31,6,1)$, $(57,8,1)$, and $(73,9,1)$ difference sets were constructed,
corresponding to projective planes of orders $2$, $5$, $7$, and~$8$.
It was proved that tilings of groups with $(13,4,1)$ difference sets,
corresponding to projective planes of order $3$ do not exist. In
Section~\ref{sec4} we give an example of a mosaic of $2$-$(13,4,1)$ designs,
thereby showing that homogenous mosaics of symmetric designs are
indeed more general than tilings of groups with difference sets.

\section{An inhomogenous mosaic}\label{sec2}

Table~\ref{table1} gives an example of a mosaic with parameters
$$2\mbox{-}(13,3,1)\oplus 2\mbox{-}(13,4,2)\oplus 2\mbox{-}(13,6,5).$$
Evidently the matrix is invariant under a cyclic rotation of the $13$ rows, columns
$1$-$13$, and columns $14$-$26$. Generally, an \emph{automorphism} of a mosaic
$M=[m_{ij}]$ is a triple of permutations $(\alpha,\beta,\gamma)\in S_v \times S_b \times S_c$
such that $\gamma(m_{\alpha(i)\beta(j)})=m_{ij}$ holds for all $i\in \{1,\ldots,v\}$ and
$j \in \{1,\ldots,b\}$. Here, $\alpha=(1,2,\ldots,13)$, $\beta=(1,2,\ldots,13)(14,15,\ldots,26)$,
and $\gamma=()$ is the identity.

\begin{table}[t]
$$\left[\begin{array}{@{\,\,}c@{\,\,\,\,}c@{\,\,\,\,}c@{\,\,\,\,}c@{\,\,\,\,}c@{\,\,\,\,}c@{\,\,\,\,}c@{\,\,\,\,}
c@{\,\,\,\,}c@{\,\,\,\,}c@{\,\,\,\,}c@{\,\,\,\,}c@{\,\,\,\,}c@{\,\,\,\,}c@{\,\,\,\,}c@{\,\,\,\,}c@{\,\,\,\,}
c@{\,\,\,\,}c@{\,\,\,\,}c@{\,\,\,\,}c@{\,\,\,\,}c@{\,\,\,\,}c@{\,\,\,\,}c@{\,\,\,\,}c@{\,\,\,\,}c@{\,\,\,\,}c@{\,\,}}
1 & 3 & 3 & 3 & 2 & 3 & 2 & 2 & 3 & 1 & 3 & 2 & 1 & 1 & 3 & 2 & 2 & 3 & 3 & 1 & 3 & 3 & 3 & 2 & 1 & 2\\
1 & 1 & 3 & 3 & 3 & 2 & 3 & 2 & 2 & 3 & 1 & 3 & 2 & 2 & 1 & 3 & 2 & 2 & 3 & 3 & 1 & 3 & 3 & 3 & 2 & 1\\
2 & 1 & 1 & 3 & 3 & 3 & 2 & 3 & 2 & 2 & 3 & 1 & 3 & 1 & 2 & 1 & 3 & 2 & 2 & 3 & 3 & 1 & 3 & 3 & 3 & 2\\
3 & 2 & 1 & 1 & 3 & 3 & 3 & 2 & 3 & 2 & 2 & 3 & 1 & 2 & 1 & 2 & 1 & 3 & 2 & 2 & 3 & 3 & 1 & 3 & 3 & 3\\
1 & 3 & 2 & 1 & 1 & 3 & 3 & 3 & 2 & 3 & 2 & 2 & 3 & 3 & 2 & 1 & 2 & 1 & 3 & 2 & 2 & 3 & 3 & 1 & 3 & 3\\
3 & 1 & 3 & 2 & 1 & 1 & 3 & 3 & 3 & 2 & 3 & 2 & 2 & 3 & 3 & 2 & 1 & 2 & 1 & 3 & 2 & 2 & 3 & 3 & 1 & 3\\
2 & 3 & 1 & 3 & 2 & 1 & 1 & 3 & 3 & 3 & 2 & 3 & 2 & 3 & 3 & 3 & 2 & 1 & 2 & 1 & 3 & 2 & 2 & 3 & 3 & 1\\
2 & 2 & 3 & 1 & 3 & 2 & 1 & 1 & 3 & 3 & 3 & 2 & 3 & 1 & 3 & 3 & 3 & 2 & 1 & 2 & 1 & 3 & 2 & 2 & 3 & 3\\
3 & 2 & 2 & 3 & 1 & 3 & 2 & 1 & 1 & 3 & 3 & 3 & 2 & 3 & 1 & 3 & 3 & 3 & 2 & 1 & 2 & 1 & 3 & 2 & 2 & 3\\
2 & 3 & 2 & 2 & 3 & 1 & 3 & 2 & 1 & 1 & 3 & 3 & 3 & 3 & 3 & 1 & 3 & 3 & 3 & 2 & 1 & 2 & 1 & 3 & 2 & 2\\
3 & 2 & 3 & 2 & 2 & 3 & 1 & 3 & 2 & 1 & 1 & 3 & 3 & 2 & 3 & 3 & 1 & 3 & 3 & 3 & 2 & 1 & 2 & 1 & 3 & 2\\
3 & 3 & 2 & 3 & 2 & 2 & 3 & 1 & 3 & 2 & 1 & 1 & 3 & 2 & 2 & 3 & 3 & 1 & 3 & 3 & 3 & 2 & 1 & 2 & 1 & 3\\
3 & 3 & 3 & 2 & 3 & 2 & 2 & 3 & 1 & 3 & 2 & 1 & 1 & 3 & 2 & 2 & 3 & 3 & 1 & 3 & 3 & 3 & 2 & 1 & 2 & 1\\
\end{array}\right]$$
\vskip 3mm
\caption{A $2\mbox{-}(13,3,1)\oplus 2\mbox{-}(13,4,2)\oplus 2\mbox{-}(13,6,5)$ mosaic.}\label{table1}
\end{table}

The designs in this mosaic are obtained by developing the ordered difference families $\F_1=(\{0, 1, 4 \},\,
\{0, 2, 7\})$, $\F_2=(\{2, 6, 7, 9\},\, \{ 1, 3, 10, 11 \})$, and $\F_3 = (\{ 3, 5, 8, 10, 11, 12 \},\,
\{ 4, 5, 6, 8, 9, 12 \})$ in $\Z_{13}$. The first components of $\F_1$, $\F_2$, and $\F_3$ consitute
a partition of $\Z_{13}$, as do the second components.

\section{Homogenous mosaics of designs that are not resolvable}\label{sec3}

Table~\ref{table2} gives two homogenous mosaics of $2$-$(9,3,2)$ designs. According to
Table~24 in~\cite[Section~II.1.3]{CD07}, there are $36$ designs with these parameters
up to isomorphism, $9$ of them resolvable. The first mosaic contains three isomorphic
copies of a non-resolvable $2$-$(9,3,2)$ design. The second mosaic contains three non-isomorphic
design. The design represented by entries~$1$ is resolvable, while the designs represented
by entries~$2$ and~$3$ are not resolvable. Hence, these two mosaics cannot be obtained
by \cite[Theorem~3.4]{GGP18}.

\begin{table}[!ht]
$$\left[\begin{array}{@{\,\,}c@{\,\,\,\,}c@{\,\,\,\,}c@{\,\,\,\,}c@{\,\,\,\,}c@{\,\,\,\,}c@{\,\,\,\,}c@{\,\,\,\,}
c@{\,\,\,\,}c@{\,\,\,\,}c@{\,\,\,\,}c@{\,\,\,\,}c@{\,\,\,\,}c@{\,\,\,\,}c@{\,\,\,\,}c@{\,\,\,\,}c@{\,\,\,\,}
c@{\,\,\,\,}c@{\,\,\,\,}c@{\,\,\,\,}c@{\,\,\,\,}c@{\,\,\,\,}c@{\,\,\,\,}c@{\,\,\,\,}c@{\,\,}}
1 & 2 & 1 & 1 & 2 & 1 & 1 & 3 & 3 & 1 & 2 & 3 & 1 & 3 & 2 & 1 & 3 & 3 & 2 & 2 & 3 & 2 & 3 & 2 \\
1 & 1 & 2 & 1 & 1 & 2 & 3 & 1 & 3 & 3 & 1 & 2 & 2 & 1 & 3 & 3 & 1 & 3 & 3 & 2 & 2 & 2 & 2 & 3 \\
2 & 1 & 1 & 2 & 1 & 1 & 3 & 3 & 1 & 2 & 3 & 1 & 3 & 2 & 1 & 3 & 3 & 1 & 2 & 3 & 2 & 3 & 2 & 2 \\
1 & 3 & 2 & 2 & 3 & 3 & 1 & 2 & 1 & 3 & 3 & 1 & 2 & 1 & 2 & 2 & 2 & 3 & 1 & 3 & 1 & 1 & 3 & 2 \\
2 & 1 & 3 & 3 & 2 & 3 & 1 & 1 & 2 & 1 & 3 & 3 & 2 & 2 & 1 & 3 & 2 & 2 & 1 & 1 & 3 & 2 & 1 & 3 \\
3 & 2 & 1 & 3 & 3 & 2 & 2 & 1 & 1 & 3 & 1 & 3 & 1 & 2 & 2 & 2 & 3 & 2 & 3 & 1 & 1 & 3 & 2 & 1 \\
2 & 3 & 3 & 1 & 3 & 2 & 3 & 2 & 2 & 2 & 2 & 1 & 1 & 3 & 3 & 2 & 1 & 1 & 1 & 2 & 3 & 3 & 1 & 1 \\
3 & 2 & 3 & 2 & 1 & 3 & 2 & 3 & 2 & 1 & 2 & 2 & 3 & 1 & 3 & 1 & 2 & 1 & 3 & 1 & 2 & 1 & 3 & 1 \\
3 & 3 & 2 & 3 & 2 & 1 & 2 & 2 & 3 & 2 & 1 & 2 & 3 & 3 & 1 & 1 & 1 & 2 & 2 & 3 & 1 & 1 & 1 & 3 \\
\end{array}\right]$$

$$\left[\begin{array}{@{\,\,}c@{\,\,\,\,}c@{\,\,\,\,}c@{\,\,\,\,}c@{\,\,\,\,}c@{\,\,\,\,}c@{\,\,\,\,}c@{\,\,\,\,}
c@{\,\,\,\,}c@{\,\,\,\,}c@{\,\,\,\,}c@{\,\,\,\,}c@{\,\,\,\,}c@{\,\,\,\,}c@{\,\,\,\,}c@{\,\,\,\,}c@{\,\,\,\,}
c@{\,\,\,\,}c@{\,\,\,\,}c@{\,\,\,\,}c@{\,\,\,\,}c@{\,\,\,\,}c@{\,\,\,\,}c@{\,\,\,\,}c@{\,\,}}
1 & 2 & 1 & 1 & 2 & 1 & 1 & 3 & 3 & 1 & 3 & 3 & 1 & 3 & 2 & 1 & 2 & 3 & 3 & 2 & 2 & 3 & 2 & 2 \\
1 & 1 & 2 & 1 & 1 & 2 & 3 & 1 & 3 & 3 & 1 & 3 & 2 & 1 & 3 & 3 & 1 & 2 & 2 & 3 & 2 & 2 & 3 & 2 \\
2 & 1 & 1 & 2 & 1 & 1 & 3 & 3 & 1 & 3 & 3 & 1 & 3 & 2 & 1 & 2 & 3 & 1 & 2 & 2 & 3 & 2 & 2 & 3 \\
1 & 3 & 2 & 3 & 3 & 1 & 2 & 2 & 1 & 2 & 1 & 3 & 3 & 3 & 2 & 2 & 3 & 2 & 1 & 2 & 1 & 1 & 3 & 1 \\
2 & 1 & 3 & 1 & 3 & 3 & 1 & 2 & 2 & 3 & 2 & 1 & 2 & 3 & 3 & 2 & 2 & 3 & 1 & 1 & 2 & 1 & 1 & 3 \\
3 & 2 & 1 & 3 & 1 & 3 & 2 & 1 & 2 & 1 & 3 & 2 & 3 & 2 & 3 & 3 & 2 & 2 & 2 & 1 & 1 & 3 & 1 & 1 \\
2 & 3 & 3 & 2 & 3 & 2 & 2 & 3 & 1 & 1 & 2 & 2 & 1 & 1 & 2 & 3 & 1 & 1 & 1 & 3 & 3 & 3 & 1 & 2 \\
3 & 2 & 3 & 2 & 2 & 3 & 1 & 2 & 3 & 2 & 1 & 2 & 2 & 1 & 1 & 1 & 3 & 1 & 3 & 1 & 3 & 2 & 3 & 1 \\
3 & 3 & 2 & 3 & 2 & 2 & 3 & 1 & 2 & 2 & 2 & 1 & 1 & 2 & 1 & 1 & 1 & 3 & 3 & 3 & 1 & 1 & 2 & 3 \\
\end{array}\right]$$
\vskip 3mm
\caption{Two $2\mbox{-}(9,3,2)\oplus 2\mbox{-}(9,3,2)\oplus 2\mbox{-}(9,3,2)$ mosaics.}\label{table2}
\end{table}

We constructed these mosaics by assuming an automorphism $\alpha=(1,2,3)$ $(4,5,6)(7,8,9)$,
$\beta=(1,2,3)(4,5,6)\cdots (22,23,24)$, $\gamma=()$ of order~$3$ and using a modified
version of the Kramer-Mesner method~\cite{KM76}. We found more homogenous mosaics like this
with parameters $2$-$(12,4,3)$ and automorphisms of order~$11$.

\section{A mosaic of projective planes of order $3$}\label{sec4}

Perhaps our most interesting example is given in Table~\ref{table3}.
This is a homogenous mosaic of symmetric $2$-$(13,4,1)$ designs,
i.e.\ projective planes of order~$3$. We did an exhaustive computer
search and found only this example up to isomorphism and transposition.
Its full automorphism group is of order~$3$, generated by
$\alpha=\beta=(1,2,3)(4,5,6)(7,8,9)(10,11,12)$ and $\gamma=(1,2,3)$.
Mosaics obtained from tilings of a group~$G$ by difference sets have~$G$
as an automorphism group acting regularly on the rows and columns,
i.e.\ points and blocks of the designs. Hence, this mosaic
does not come from a tiling of the group~$\Z_{13}$ by $(13,4,1)$
difference sets.

\begin{table}[!ht]
$$\left[\begin{array}{@{\,\,}c@{\,\,\,\,}c@{\,\,\,\,}c@{\,\,\,\,}c@{\,\,\,\,}c@{\,\,\,\,}c@{\,\,\,\,}c@{\,\,\,\,}
c@{\,\,\,\,}c@{\,\,\,\,}c@{\,\,\,\,}c@{\,\,\,\,}c@{\,\,\,\,}c@{\,\,}}
0 & 1 & 2 & 1 & 3 & 2 & 3 & 1 & 1 & 3 & 3 & 2 & 2 \\
3 & 0 & 2 & 3 & 2 & 1 & 2 & 1 & 2 & 3 & 1 & 1 & 3 \\
3 & 1 & 0 & 2 & 1 & 3 & 3 & 3 & 2 & 2 & 1 & 2 & 1 \\
3 & 3 & 1 & 0 & 1 & 1 & 2 & 2 & 1 & 2 & 3 & 3 & 2 \\
2 & 1 & 1 & 2 & 0 & 2 & 2 & 3 & 3 & 1 & 3 & 1 & 3 \\
2 & 3 & 2 & 3 & 3 & 0 & 1 & 3 & 1 & 2 & 2 & 1 & 1 \\
1 & 2 & 2 & 2 & 3 & 3 & 0 & 2 & 1 & 1 & 1 & 3 & 3 \\
3 & 2 & 3 & 1 & 3 & 1 & 2 & 0 & 3 & 1 & 2 & 2 & 1 \\
1 & 1 & 3 & 2 & 2 & 1 & 1 & 3 & 0 & 3 & 2 & 3 & 2 \\
1 & 3 & 3 & 1 & 1 & 2 & 3 & 2 & 2 & 0 & 2 & 1 & 3 \\
1 & 2 & 1 & 3 & 2 & 2 & 3 & 1 & 3 & 2 & 0 & 3 & 1 \\
2 & 2 & 3 & 3 & 1 & 3 & 1 & 1 & 2 & 1 & 3 & 0 & 2 \\
2 & 3 & 1 & 1 & 2 & 3 & 1 & 2 & 3 & 3 & 1 & 2 & 0 \\
\end{array}\right]$$
\vskip 3mm
\caption{A $2\mbox{-}(13,4,1)\oplus 2\mbox{-}(13,4,1)\oplus 2\mbox{-}(13,4,1)\oplus 2\mbox{-}(13,1,0)$ mosaic.}\label{table3}
\end{table}

In~\cite{CKZ15} it was proved that tilings of groups with difference sets are
not possible for $(13,4,1)$, $(15,7,3)$, $(21,5,1)$, $(35,17,8)$, and $(40,13,4)$.
Of course, they are also not possible if there are no difference sets, as is
the case for $(31,10,3)$. Symmetric designs with all of these parameters
exist and it might be possible to constructs homogenous mosaics,
as we did for $(13,4,1)$.

The examples presented in this paper are available on the web page
\begin{center}
\url{https://web.math.pmf.unizg.hr/~krcko/results/mosaics.html}
\end{center}
They can be analyzed using the GAP package \emph{Prescribed Automorphism
Groups}~\cite{PAG} and our claims about their properties can be verified.

\end{document}